\documentclass[12pt]{amsart}

\newcommand {\supplus}{\mathop{{\supset}\llap{\raise
0.5pt\hbox{\normalfont\small+}\hskip 0.5pt}}}

\newcommand {\subplus}{\mathop{{\subset}\llap{\raise
0.5pt\hbox{\normalfont\small+}\hskip 0.5pt}}}

\newcommand {\Cee}    {{\mathbb  C}}

\newcommand {\fA}     {{\mathfrak{A}}}

\newcommand {\fb}     {{\mathfrak{b}}}

\newcommand {\fg}     {{\mathfrak{g}}}    %
\newcommand {\fgl}    {{\mathfrak{gl}}}  %
\newcommand {\fh}     {{\mathfrak{h}}}

\newcommand {\fq}     {{\mathfrak{q}}}

\newcommand {\fS}     {{\mathfrak{S}}}

\newcommand {\cal} {\mathcal}

\def \opname#1#2%
  {\expandafter\newcommand \csname #1\endcsname {{\mathop{#2}\nolimits}}}

\newcommand{\rmname}[1]
  {\expandafter\newcommand \csname #1\endcsname {{\operatorname{#1}}}}

\newcommand{\rmnameii}[2]
  {\expandafter\newcommand \csname #1\endcsname {{\operatorname{#2}}}}

\rmname{act}
\rmname{Ad}
\rmname{Add}
\rmname{ad}
\rmname{Alt}
\rmname{alt}
\rmname{Ann}
\rmname{antidiag}
\rmname{Ber}
\rmname{ber}
\rmname{Br}
\rmname{card}
\rmname{ch}
\rmname{Char}
\rmname{cem}
\rmname{cj}
\rmname{Cliff}
\rmname{cntr}
\rmname{codim}
\rmname{coind}
\rmname{const}
\rmname{col}
\rmname{cork}
\rmname{cpr}
\rmname{diag}
\rmnameii{Div}{div}
\rmname{Def}
\rmname{Der}
\rmname{Dim}
\rmname{End}
\rmname{Even}
\rmname{Ext}
\rmname{gr}
\rmname{Hom}
\rmname{HT}
\rmnameii{Ht}{ht}
\rmname{hwt}
\rmname{Id}
\rmname{id}
\rmname{ind}
\rmname{Ind}
\rmname{Inf}
\rmname{irr}
\rmname{Le}
\rmname{Lie}
\rmname{lwt}
\rmname{mult}
\rmname{Mor}
\rmname{nm}
\rmname{Ob}
\rmname{Odd}
\rmname{Osc}
\rmname{per}
\rmname{Pic}
\rmname{pr}
\rmname{pro}
\rmname{Prime}
\rmname{Proj}
\rmname{prt}
\rmname{pt}
\rmname{Q}
\rmname{qet}
\rmname{qtr}
\rmname{rd}
\rmname{rk}
\rmname{row}
\rmname{Res}
\rmname{salt}
\rmname{Sch}
\rmname{SBr}
\rmname{scalar}
\rmname{Ser}
\rmname{sign}
\rmname{Smbl}
\rmname{spin}
\rmname{ssym}
\rmname{str}
\rmname{st}
\rmname{sgn}
\rmname{sq}
\rmname{symm}
\rmname{supp}
\rmname{Supp}
\rmname{St}
\rmname{Spec}
\rmname{Spm}
\rmname{tr}
\rmname{vpt}
\rmname{weyl}
\rmname{Weyl}
\rmname{Witt}

\opname{vvol}  {{v\hspace{-0.1ex}o\hspace{-0.02ex}l\/}}
\opname{pnt}  {\text{\normalfont pt}}
\opname{Span} {{Span}}
\opname{slim} {\overline{\lim}}
\opname{Vol}  {{V\hspace{-0.55ex}o\hspace{-0.02ex}l\/}}
\opname{Par}  {{P\hspace{-0.3ex}a\hspace{-0.05ex}r\/}}

\rmname{Mat}
\rmname{Bil}
\rmname{Diff}
\rmname{Ker}
\rmname{Herm}
\rmname{Coker}
\rmname{Conn}
\rmname{Covect}
\rmname{Vect}
\rmname{Int}

\rmnameii {IM} {Im}
\rmnameii {RE} {Re}

\opname{Aut} {{A\hspace{-0.2ex}u\hspace{-0.1ex}t\/}}
\opname{GL} {{G\hspace{-0.3ex}L}}
\opname{SL} {{S\hspace{-0.3ex}L}}
\opname{Exp} {{E\hspace{-0.2ex}x\hspace{-0.1ex}p\/}}
\opname{GQ} {{G\hspace{-0.2ex}Q}}
\opname{OSp} {{O\hspace{-0.25ex}S\hspace{-0.15ex}p\/}}
\opname{Out} {{O\hspace{-0.25ex}u\hspace{-0.15ex}t\/}}
\opname{Spp} {{S\hspace{-0.2ex}p\/}}
\opname{SpO} {{S\hspace{-0.2ex}p\hspace{-0.02ex}O\/}}
\opname{Pe} {{P\hspace{-0.25ex}e\/}}
\opname{SPe} {{S\hspace{-0.25ex}P\hspace{-0.25ex}e\/}}
\opname{Spin} {{S\hspace{-0.25ex}p\hspace{-0.05ex}i\hspace{-0.1ex}n\/}}
\opname{Iso} {{I\hspace{-0.25ex}s\hspace{-0.1ex}o\/}}
\opname{SSPe} {{S\hspace{-0.25ex}S\hspace{-0.15ex}P\hspace{-0.25ex}e\/}}
\opname{PeU} {{P\hspace{-0.25ex}e\hspace{-0.1ex}U\/}}
\opname{QU} {{Q\hspace{-0.15ex}U\/}}
\opname{U} {{U\/}}

\opname{cGQ} {{\cal G \hspace{-0.2em} Q \/}}
\opname{cSL} {{\cal S \hspace{-0.2em} L \/}}
\opname{cGL} {{\cal G \hspace{-0.2em} L \/}}
\opname{cOSp} {{\cal O \hspace{-0.2em} S \hspace{-0.3em} \it p\/}}
\opname{cPe} {{\cal P \hspace{-1.5pt} \it e\/}}
\opname{cVect} {{\cal V \hspace{-1.5pt} \it
e\hspace{-0.1ex}c\hspace{-0.1ex}t\/}}
\opname{cVol} {{\cal V \hspace{-1.5pt} \it o\hspace{-0.1ex}l\/}}
\opname{cAut} {{\cal A \hspace{-0.2em} \it u\hspace{-0.1em}t\/}}
\opname{cCovect} {{\cal C \hspace{-1.5pt}
     \it o\hspace{-0.1ex}v\hspace{-0.1ex}e\hspace{-0.1ex}c\hspace{-0.1ex}t\/}}
\opname{CW} {{C\hspace{-0.15ex}W}}

\newcommand {\ev} {{\bar0}}
\newcommand {\od} {{\bar1}}

\newcommand {\tto} {\longrightarrow}

\newcommand {\bcdot}   {\mathbin{\hbox{\raise.4ex\hbox{\bf.}}}} 

\newcommand {\secno} {}
\newcommand {\ssecfont} {\normalfont\bf}

\newtheorem{Theorem}{\secno Theorem}

\newenvironment {th*}[1]
    {\gdef\thname{#1} \begin{thn}}%
    {\end{thn}}
\newtheorem{thn}[Theorem] {\thname}

\theoremstyle{definition}

\newenvironment {ex*}[1]
    {\gdef\thname{#1} \begin{exn}}%
    {\end{exn}}
\newtheorem{exn}[Theorem]{\thname}

\theoremstyle{remark}

\newenvironment {rem*}[1]
    {\gdef\thname{#1} \begin{remn}}%
    {\end{remn}}
\newtheorem{remn}[Theorem]{\thname}

\newcommand {\ssec}{\subsection*}

\newcommand {\ssbegin}[2]
  {\def \secno {\gdef \secno {}{\ssecfont #1. }}%
   \begin{#2}}
\begin{document}

\title[Projective Representations of Symmetric Groups]{The Howe
duality and \\ the Projective Representations of Symmetric Groups}

\author{Alexander Sergeev}

\address{Dept.  of Math., Univ.  of Stockholm, Roslagsv.  101,
Kr\"aftriket hus 6, S-106 91, Stockholm, Sweden. On leave of absence
from Balakovo Inst.  of Technique Technology and Control\\ e-mail:
mleites@matematik.su.se (subject: for Sergeev)}

\thanks{I am thankful to D. Leites for support and help.}

\begin{abstract} The symmetric group $\fS_{n}$ possesses a nontrivial
central extension, whose irreducible representations, different from
the irreducible representations of $\fS_{n}$ itself, coincide with the
irreducible representations of the algebra $\fA_{n}$ generatored by
indeterminates $\tau_{i, j}$ for $i\neq j$, $1\leq i, j\leq n$ subject to
the relations
$$
\renewcommand{\arraystretch}{1.4}
\begin{array}{l}
	\tau_{i, j}=-\tau_{j, i}, \quad \tau_{i, j}^{2}=1, \quad \tau_{i,
j}\tau_{k, l}=-\tau_{k, l}\tau_{i, j}\text{ if }\{i, j\}\cap\{k,
l\}=\emptyset;\\
\tau_{i, j}\tau_{j, k}\tau_{i, j}=\tau_{j, k}\tau_{i, j}\tau_{j,
k}=-\tau_{i, k}\; \; \text{ for  any } i, j, k, l.\\
\end{array}
$$
Recently M.~Nazarov realized irreducible representations of $\fA_{n}$
and Young symmetrizers by means of the Howe duality between the Lie
superalgebra $\fq(n)$ and the Hecke algebra $H_{n}=\fS_{n}\circ
C_{n}$, the semidirect product of $\fS_{n}$ with the Clifford algebra
$C_{n}$ on $n$ indeterminates.

Here I construct one more analog of Young symmetrizers in $H_{n}$ as
well as the analogs of Specht modules for $\fA_{n}$ and $H_{n}$.
\end{abstract}



\maketitle

\section*{\S 1.  Introduction} Lately we witness an increase of
interest in the study of representations of symmetric groups.  In
particular, in their projective representations.

Recall that the symmetric group $\fS_{n}$ has a nontrivial central
extension whose irreducible representations do not reduce to those of
$\fS_{n}$ but coincide (may be identified) with the irreducible
representations of the algebra $\fA_{n}$ determined by generators
$\tau_{i, j}$ for $i\neq j$, $1\leq i, j\leq n$ subject to the
relations
$$
\renewcommand{\arraystretch}{1.4}
\begin{array}{l}
	\tau_{i, j}=-\tau_{j, i}, \quad \tau_{i, j}^{2}=1, \quad \tau_{i,
j}\tau_{k, l}=-\tau_{k, l}\tau_{i, j}\text{ if }\{i, j\}\cap\{k,
l\}=\emptyset;\\
\tau_{i, j}\tau_{j, k}\tau_{i, j}=\tau_{j, k}\tau_{i, j}\tau_{j,
k}=-\tau_{i, k}\; \; \text{ for  any } i, j, k, l.\\
\end{array}\eqno{(1.1)}
$$
In \cite{N1} Nazarov realized irreducible representations of $\fA_{n}$
by means of an orthogonal basis constructed in each of the spaces of
the representations and indicating the action of the generators
$\tau_{i, i+1}$ on them (an analog of the Young orthogonal form).  In
\cite{N2}, with the help of an ``odd'' analog of the degenerate affine
Hecke algebra Nazarov constructed elements of the algebra
$H_{n}=\fS_{n}\circ  C_{n}$, the semidirect product of $\fS_{n}$ with
the Clifford algebra $C_{n}$ on $n$ indeterminates.  The elements of
$H_{n}$ serve as analogs of Young symmetrizers.

Here I construct one more analog of Young symmetrizers in $H_{n}$ as
well as analogs of Specht modules (cf.  \cite{Ja}) for the algebras
$\fA_{n}$ and $H_{n}$.  This construction is based on another form of
expression of Young symmetrizers for $\fS_{n}$.

Namely, let $t$ be a Young tableau
(i.e., a Young diagram filled in with numbers 1 to $n$), $R_{t}$ and
$C_{t}$ the row and column stabilizers of $t$;
$$
\rho_{t}=\sum_{\sigma\in R_{t}}\sigma,\quad \kappa_{t}=
\sum_{\sigma\in C_{t}}\varepsilon(\sigma)\sigma.
$$
Then, up to a constant factor, the Young symmetrizer of $t$ is
$$
e_{t}= \kappa_{t}\rho_{t}.\eqno{(1.2)}
$$

Let us represent the Young symmetrizer differently.  Let $I$ be the
sequence obtained by reading the tableau $t$ along columns left to
right and downwards.  For each $i\in I$ the {\it Jucys--Murphy
elements} (\cite{J}, \cite{M}) are defined to be
$$
x_{i}=\sum_{\alpha\text{ preceeds $i$}}s_{\alpha i}, \text{
where the $s_{\alpha i}\in\fS_{n}$ are transpositions}.
$$

One can verify that the $x_{i}$ commute with each other. Set
$$
\tilde \kappa_{t}=\prod_{i\in I}(j-x_{i}), \;\text{ where $j$ is
the number of the column occupied by $i$}.
$$
It is subject to a direct verification that $e_{t}$ can be also
expressed as
$$
e_{t}= \tilde \kappa_{t}\rho_{t}.\eqno{(1.3)}
$$

I use this other representation of the Young symmetrizers to construct
the corresponding elements in $H_{n}$ and prove that theyse
corresponding elements are idempotents that generate isotypical
ideals.  With the help of these elements and their analogs for
$\fA_{n}$ I construct a realization of multiples of irreducible
modules similar to Specht modules, cf.  \cite{Ja}.

The proofs are based on the Howe duality between Lie superalgebra
$\fq(n)$ and $H_{n}$.

\section*{\S 2. Background}
Let $\fS_{n}$ be the symmetric group, $C_{n}$ the Clifford algebra
generated by $n$ indeterminates $p_{1}$, \dots , $p_{n}$ subject to
the relations
$$
p_{i}^2=-1,\; p_{i}p_{j}+p_{j}p_{i}=0\text{ for }i\neq j.
$$
The symmetric group acts on $C_{n}$ permuting the generators, so we
can form a semidirect product $H_{n}=\fS_{n}\circ  C_{n}$. Set
$$
\tau_{i,j}=\frac{1}{\sqrt{2}}(p_{i}-p_{j})s_{i,j}.
$$
As is not difficult to verify, the relations (1.1) hold; hence, the
algebra generated by the $\tau_{i,j}$ is isomorphic to $\fA_{n}$.
Besides, $\fA_{n}$ supercommutes with $C_{n}$; hence,
$H_{n}=\fA_{n}\otimes C_{n}$, as {\it super}algebras if we define
parity in $\fA_{n}$  by setting
$p(\tau_{i,j})=\od$.

Let $V$ be a superspace of superdimension $(n, n)$ with the fixed
basis $\{e_{i}\}_{i=1}^n\cup\{e_{\bar i}\}_{\bar i=\bar 1}^{\bar n}$ and
the odd operator
$$
Q: e_{i}\mapsto (-1)^{\bar i} e_{\bar i}; \quad e_{\bar i}\mapsto
 e_{i}.
$$
The (super)centrilizer of $Q$ in $\Mat(V)$ is denoted by $\Q(V)$, cf.
\cite{BL}.
We denote the Lie superalgebras associated with the associative
superalgebras $\Mat(V)$ and $\Q(V)$ by $\fgl(V)$ and $\fq(V)$,
respectively.

As is shown in \cite{S1}, the algebras $\fq(V)$ and $H_{k}$ constitute a
Howe-dual pair in the superspace $W=V^{\otimes k}$. Therefore, the
following decomposition takes place:
$$
W=\mathop{\oplus}\limits_{\lambda}2^{-\delta(\lambda)}T^{\lambda}
\otimes V^{\lambda},\eqno{(2.1)}
$$
where $\lambda$ runs over strict partitions of $k$, $T^{\lambda}$ is an
irreducible (in supersence) $H_{k}$-module, $V^{\lambda}$ an
irrreducible $\fq(n)$-module and $\lambda_{n+1}=0$ and where
$$
\delta(\lambda)=\left\{\begin{matrix}0&\text{if the number of
nonzero parts of $\lambda$ is even}\cr
1&\text{otherwise.}\end{matrix}\right .
$$

Select a basis in $\fq(n)$: let the $e_{i}^*$ be the left dual
basis to the $\{e_{i}\}$, $i=1, \dots n$; $\bar 1\dots \bar n$; set
$$
E_{i, j}=e_{i}\otimes e_{j}^*+ e_{\bar i}\otimes e_{\bar j}^*; \;
F_{i, j}=e_{i}\otimes e_{\bar j}^*+ e_{\bar i}\otimes e_{j}^*.
$$
Then $\fh= \Span (E_{i, i}\text{ and }F_{i, i}: i=1, \dots, n)$ is a
Cartan subalgebra in $\fq(n)$ and $\fb =\Span (E_{i, j}\text{ and }F_{i, j}:
i\leq j)$ is a Borel subalgebra, cf.  \cite{Pe}.

If $\lambda$ is a strict partition and $\lambda_{n+1}=0$, then
$\lambda$ can be interpreted as a linear functional on $\fh_{\ev}$:
set
$$
\lambda(E_{ii})=\lambda_{i}.\eqno{(2.2)}
$$
Let $R^\lambda$ be the $H_{k}$-module equal to the direct sum of
$2^{(l(\lambda)-\delta(\lambda))/2}$ copies of $T^\lambda$. It is not
difficult to see that $R^\lambda$ coincides with the set of
$\fb$-highest weight vectors of weight $\lambda$ in $W$.

\ssbegin{2.1}{Lemma} Let $\mu$ be a strict partition, $\mu_{n+1}=0$.
Let $V^{\mu}$ be an isotypical modul of type $\mu$ and $n_{\mu}$ the
multiplicity of  the highsest weight vector in $V^{\mu}$. Then the
highest weights of $V^{\mu}\otimes V$ are of the form
$\mu+\varepsilon_{i}$, where the $\varepsilon_{i}$ are the weights
of $V$,  and their multiplicity is equal to $2n_{\mu}$.
\end{Lemma}

Proof follows easily from the multiplication table of the projective
Schur functions, see \cite{P}. \qed

\ssbegin{2.2}{Lemma} Let $\fg=\fq(n)$; $\fb$ and $\fh$ be defined as
above and $V$ a $\fg$-module. Let $V^{+}_{\lambda}$ be the set of
$\fb$-highest vector of weight $\lambda$.

If $u\in U(\fg)$ and $uV^{+}_{\lambda}\subset V^{+}_{\lambda}$, then
there exists $w\in U(\fh)$ such that
$u|_{V^{+}_{\lambda}}=w|_{V^{+}_{\lambda}}$.
\end{Lemma}

\begin{proof} Let $u=\sum u_{\alpha}$ be the weight decomposition of
$u\in U(\fg)$ with respect to $\fh_{\ev}$, the even part of the Cartan
subalgebra $\fh$.  Therefore, if $v\in V_{\lambda}$, then $uv=\sum
_{\alpha} u_{\alpha}v$ and, if $u_{\alpha}v\neq 0$, then the weight of
$u_{\alpha}v$ is equal to $\lambda +\alpha$.  Thus, we may assume that
$u=u_{0}\in U(\fg)^{\fh_{\ev}}$, where $U(\fg)^{\fh_{\ev}}$ is the
centralizer of $\fh$.

Thanks to \cite{S2}, we know that $U(\fg)^{\fh_{\ev}}\cong
U(\fh)\oplus L$, where $L=U(\fg)^{\fh_{\ev}}\cap U(\fg)\fb^+$, where
$\fb^+$ is the linear span of the positive roots in $\fb$, is a
twosided ideal in $U(\fg)^{\fh_{\ev}}$.

Hence, $u=w+u_{1}$, where $w\in U(\fh)$ and $u_{1}\in L$; this
implies that $uv=wv$ for $v\in V_{\lambda}$.
\end{proof}

\section*{\S 3. Specht modules over $H_{k}$}
Let $\lambda$ be a strict partition and $t$ the shifted tableau of the
form $\lambda$, where $\lambda_{1}>\lambda_{2}>\dots >\lambda_{n}>0$
and $\sum \lambda_{i}=k$. Let us fill in the tableau with the numbers
1 to $k$ and define the functional $\lambda$ on the Cartan
subalgebra $\fh$ by setting
$$
\lambda(E_{ii})=\lambda_{i}; \quad \lambda(F_{ii})=0.
$$
In $W=V^{\otimes k}$, where $V$ is the standard $\fq(V)$-module of
dimension $(n, n)$, consider the submodule $M^{\lambda}$ consisting
of the vectors of weight $\lambda$.

Again, let $R_{t}$ be the row stabilizer of $t$ and
$\rho_{t}=\sum_{\sigma\in R_{t}}\sigma$. Let $I$ be the sequence
obtained by reading the tableau $t$ downwards and from left to right.
For $i\in I$ define: $\pi_{i}=\sum \tau_{\alpha, i}$, where the sum
runs over all the $\alpha$'s, $\alpha\in I$,  that preceed $i\in I$.

It is subject to a direct verification that
$$
\pi_{i}\pi_{j}+\pi_{j}\pi_{i}=0\text{ for }i\neq j.
$$
These are odd analogs of the Jucys-Murphy elements. Set
$$
\kappa_{t}=\prod_{i\in I}(\frac{1}{2}j(j+1)-\pi_{i}^2),
$$
where $j$ is the number of the column occupied by $i$.

\ssbegin{3.1}{Theorem} Let $v_{t}\in M^\lambda$ be the vector whose
stabilizer is $R_{t}$. Then
$\kappa_{t}(v_{t})\neq 0$, $\kappa_{t}(v_{t})\in R^\lambda$ and
$\kappa_{t}(M^\lambda)=C_{k} \kappa_{t}(v_{t})$.
\end{Theorem}

\begin{proof} Induction on $\sum \lambda_{i}=k$. It suffices to assume
that $t$ is filled in consequently columnwise, from left to right and
downwards. Let $s$ be the tableau obtained from $t$ by deleting the
last cell and $\mu$ the corresponding partition while $l$ is the
length of the last column and $r$ is its number. Then $\kappa_{t}=
\kappa_{s}\cdot \kappa_{k}$, where $\kappa_{s}$ corresponds
to tableau $s$ and $\kappa_{k}=\frac{1}{2}r(r+1)-\pi_{k}^2$.

By induction, $\kappa_{s}(M^\mu)=C_{k-1} \kappa_{s}(v_{s})\subset
R^\mu$. Therefore, $\kappa_{s}(M^\nu)=0$ for any $\nu>\mu$ ordered
with respect to {\it dominance}, cf. \cite{M}. Hence,
$$
\kappa_{s}(M^\lambda)=\mathop{\oplus}\limits_{i=1}^{n}\kappa_{s}(M^{\lambda-
\varepsilon _{i}})\otimes e_{i}.
$$
If $i>l$, then $\lambda-\varepsilon _{i}>\lambda-\varepsilon
_{l}=\mu$ and, by induction,
$\kappa_{s}(M^{\lambda-\varepsilon _{i}})=0$. Hence,
$$
\kappa_{s}(M^\lambda)=\mathop{\oplus}\limits_{i=1}^{l}\kappa_{s}(M^{\lambda-
\varepsilon _{i}})\otimes e_{i}.
$$

The inequality $\nu\geq \lambda-\varepsilon _{i}$ true for $i<l$
implies that $\nu\geq \mu$, so the same irreducible $H_{k-1}$-modules
enter the decomposition of $M^{\lambda-\varepsilon _{i}}$ as those
that enter that of $M^{\mu}$.  Therefore, if $m_{i}\in M^{\lambda-
\varepsilon _{i}}$ and $\kappa_{s}(m_{i})\neq 0$, then there exist
$m\in M^\mu$ and a homomorphism $\varphi_{i}: M^{\mu}\tto
M^{\lambda-\varepsilon _{i}}$ such that
$$
\varphi_{i}(m)=m_{i}.
$$
Applying $\kappa_{s}$ to this identity we get
$$
\kappa_{s}(m_{i})=\kappa_{s}(\varphi_{i}(m))=
\varphi_{i}(\kappa_{s}(m))\in C_{k-1} \kappa_{s}(v_{s}).
$$

The Howe duality between $\fq(n)$ and $H_{k-1}$ allows us to assume
that $\varphi_{i}\in U(\fq(n))$. This proves that
$\kappa_{s}(M^\lambda)\subset V^\mu\otimes V$, where $V^\mu$ is
the $\fq(n)$-submodule of $V^{\otimes (k-1)}$ generated by
$\kappa_{s}(M^\mu)$. By Lemma 2.1 the highest weights of $V^\mu\otimes
V$ are of the form $\mu+\varepsilon_{i}$, so the possible weights are
only
$$
\mu+\varepsilon_{1}, \;  \mu+\varepsilon_{l}=\lambda\text{ and }
\mu+\varepsilon_{i}<\mu+\varepsilon_{l} \text{ if }i>l.
$$
So the submodule generated by the weights $\mu+\varepsilon_{i}$ for
$i>l$ does not contain
weight $\lambda$; hence, neither does it contain
$\kappa_{s}(M^\lambda)$.

Therefore, $\kappa_{s}(M^\lambda)$ is contained in the submodule
generated by the highest weight vectors of weight
$\mu+\varepsilon_{1}$ and $\mu+\varepsilon_{l}$. By Lemma 2.1 all the
highest weight vectors form a free $C_{k}$-module with two generators
whose weights are $\mu+\varepsilon_{1}$ and $\mu+\varepsilon_{l}$; so in each
of these submodules the operator $\pi_{k}^2$ acts by multiplying by a
constant. Let $v\in V^{\otimes (k-1)}$ and $e_{i}\in V$. Then it is
not difficult to verify that
$$
\pi_{k}(v\otimes e_{i})=\frac{1}{\sqrt{2}}\sum_{1\leq i\leq
n}(F_{ij}-p_{k}E_{ij}(v\otimes e_{i})),
$$
where $p_{i}$ is the change of parity operator in the $i$-th factor
of $V^{\otimes k}$.

If $v$ is a highest weight vector of weight $\mu$, than $v\otimes
e_{1}$ is a highest weight vector of weight $\mu+\varepsilon_{1}$ and
$$
\pi_{k}(v\otimes e_{1})=\frac{1}{\sqrt{2}}(F_{11}-
p_{k}E_{11})(v\otimes e_{1}),
$$
hence,
$$
\renewcommand{\arraystretch}{1.4}
\begin{array}{l}
\pi_{k}^2(v\otimes e_{1})=-\frac{1}{2}(F_{11}-
p_{k}E_{11})^2(v\otimes e_{1})=\\
\frac{1}{2}(E_{11}^2-E_{11})^2(v\otimes e_{1})=\\
\frac{1}{2}((r+1)^2-(r+1))(v\otimes e_{1})=
\frac{1}{2}r(r+1)(v\otimes e_{1}).
\end{array}
$$
Therefore, it suffices to demonstrate that
$$
\kappa_{t}(v_{t})= \kappa_{t}(v_{s}\otimes e_{l})=
\kappa_{k}(\kappa_{s}(v_{s}\otimes e_{l})\neq 0.
$$
It is not difficult to verify that if $v$ is the highest weight
vector, then
$$
\pi_{k}^2(v\otimes e_{l})=-\frac{1}{2}\sum_{k\leq
i}[(E_{ik}^{(2)}v)\otimes e_{k}+
p_{k}(F_{ik}^{(2)}v)\otimes e_{k}]+
2(\sum E_{kk}v)\otimes e_{i},
$$
where
$$
E_{ik}^{(2)}= \sum_{k\leq j\leq
i}(E_{ij}E_{jk}-F_{ij}F_{jk}),\; F_{ik}^{(2)}= \sum_{k\leq j\leq
i}(E_{ij}F_{jk}-F_{ij}E_{jk}).
$$
Therefore,
$$
\pi_{k}^2(v\otimes e_{l})=\frac{1}{2}\left (E_{ll}^{(2)}-E_{ll}+2(E_{11}+\dots
+ E_{ll})\right ) v\otimes e_{k}+\dots ,
$$
where $\dots$ replace the sum of terms of the form $u\otimes e_{i}$
with $i<l$ and
$$
\renewcommand{\arraystretch}{1.4}
\begin{array}{l}
\pi_{k}^2(v\otimes e_{l})=
\frac{1}{2}\left (\lambda_{l}^{2}-\lambda_{l}+2(\lambda_{1}+\dots
+ \lambda_{l})\right )v\otimes e_{l}+\dots=\\
\frac{1}{2}\left( (r-l)(r-l-1)+2(r+\dots
+ r-(l-1)+1+r-l)\right )v\otimes e_{l}+\dots =\\
\frac{1}{2}(r^2-r+2l-2)v\otimes e_{l}+\dots
\end{array}
$$
Consequently,
$$
\kappa_{k}(v\otimes e_{l})=
\frac{1}{2}(r^2+r-r^2+r-2l+2)(v\otimes e_{l})+\dots = (r-l+1)(v\otimes
e_{l})
+\dots \neq 0.
$$
\end{proof}

\ssbegin{3.2}{Remark} By \cite{N2} (Th.  7.2) for any strict partition
$\mu$ there exist elements $\psi_{S}\in R^\mu$, where $S$
is the standard tableau of type $\mu$, such that $p\psi_{S}$
constitute a basis of $R^\mu$ while $p$ runs over $C_{n}$.  Moreover,
the action of the elements
$$
x_{k}=\sum_{i<k}(s_{ik}+s_{ik}p_{i}p_{k})
$$
is given by the formula
$$
x_{k}\psi_{S}= C_{S}(k)\psi_{\Lambda},
$$
where $C_{S}(k)=\sqrt{(j-i)(j-i+1)}$ and  $k$ occupies the $(i, j)$-th slot.

In the above notations
$$
x_{k}=\sqrt{2}p_{k}\pi_{k},
$$
hence, $x_{k}^2=2\pi_{k}$ and
$$
\renewcommand{\arraystretch}{1.4}
\begin{array}{l}
\kappa_{t}=\mathop{\prod}\limits_{i}\frac{1}{2}\left
(j(j+1)-\pi_{i}^2\right ) =\\
\mathop{\prod}\limits_{i}\frac{1}{2}\left (j(j+1)-x_{i}^2\right )=
2^{-k}\mathop{\prod}\limits _{i}\left (j(j+1)-x_{i}^2\right ).
\end{array}
$$
As earlier, we assume that the tableau is filled in left to right and
downwards along columns.

Let us show that if $\mu\geq \lambda$, then
$\kappa_{t}\psi_{S}=0$ for any standard tableau $S$ of type
$\mu$ and distinct from $t$.

Let $k$ be the number occupying the last place of the last
column in the tableau $t$ (filled in left to right and downwards). Let
$t^{*}$ be the tableau obtained from $t$ by deleting the box with $k$;
let $S^{*}$ be similarly constructed from $S$. If $\mu^{*}\geq
\lambda^{*}$ (shapes of $S^{*}$ and $t^{*}$, respectively), then the
induction hypothsis applies. If $\mu^{*}<\lambda^{*}$, then $\mu_{1}=
\lambda_{1}+1$ and $k$ occupies the last position of the first row of
$S$. Hence, if $j$ is the number of the column of $t$ occupied by $k$,
then
$$
\begin{array}{l}
\pi_{k}^2(v_{S}) =\frac{1}{2}(j+1-1)(j+1-1+1)v_{S} =\\
\frac{1}{2}j(j+1)v_{S}.
\end{array}
$$
Hence, $\left (\frac{1}{2}j(j+1)-\pi_{k}^2\right )\psi_{S}=0$ and
$\kappa_{t}\psi_{S}=0$.

This argument shows that $\kappa_{t}$ is the projection in
$M^{\Lambda}$ onto the subspace generated by $\psi_{t}$.

At the same time
$$
\begin{array}{l}
\kappa_{k}\psi_{t}=\frac{1}{2}\left (j(j+1)-(j-i)(j-i+1)\right ) \psi_{t}=\\
\left (ij-\frac{1}{2}i(i-1)\right )\psi_{t}\neq 0.
\end{array}
$$
\end{Remark}

\ssec{3.3. Specht modules}
Recall that the {\it Specht module} for a strict partition $\lambda$
is the submodule in $M^\lambda$ generated by the
vectors $\kappa_{t}v_{t}$ for all tableaux $t$, where $v_{t}$ is the
vector whose stabilizer is equal to $R_{t}$.

\ssbegin{3.3.1}{Theorem} Specht module is equal to $R^\lambda$.  It is
isotypical and its $H_{k}$-centralizer is isomorphic to the Clifford
algebra with $l(\lambda)$ generators.
\end{Theorem}

\begin{proof} Let us realize $R^\lambda$ as the set of highest weight
vectors in $V^{\otimes k}$.  By the Howe duality between $U(\fq(V))$
and $H_{k}$ the algebra of $H_{k}$-homomorphisms is generated by
$U(\fq(V))$.  But thanks to Lemma 2.2 we may assume that the algebra
of $H_{k}$-homomorphisms is generated by $U(\fh)$ for the Cartan
subalgebra $\fh$ of $\fq(V)$. If $\lambda_{i}=0$, then in our notations
for the basis of $\fh$ we have
$$
F_{ii}^2=E_{ii}=\lambda_{i}=0 \text{ in } R^\lambda
$$
and, therefore, $\ker F_{ii}$ is an $\fh$-submodule in the set
$(V^{\lambda})^+$ of highest weight vectors. But $(V^{\lambda})^+$ is
irreducible as $\fh$-module; hence, $F_{ii}|_{(V^{\lambda})^+}=0$ which
implies $F_{ii}|_{R^{\lambda}}=0$. Thus, the algebra of $H_{k}$-homomorphisms
is generated by the $F_{ii}$ for $i$ such that $\lambda_{i}\neq 0$.

By Theorem 3.1 the Specht module is contained in $R^\lambda$ and any
homomorphism of $M^\lambda$  sends $R^\lambda$ into itself. Hence,
the Specht module coinsides with $R^\lambda$.
\end{proof}

\ssbegin{3.3.2}{Corollary} Set $p_{t}^i=\sum p_{\alpha}$, where the
$\alpha$ belongs to the $i$-th column of tableau $t$. Then for any
$H_{k}$-module endomorphism $\varphi$ of $M^\lambda$ we have
$$
\varphi(\kappa_{t}v_{t})=f\cdot \kappa_{t}v_{t},
$$
where $f$ belongs to the subalgebra of $C_{k}$ generated by
$p_{t}^1$, \dots , $p_{t}^{l(\lambda)}$.
\end{Corollary}

\begin{proof} Let us realize $M^\lambda$ as the subset of vectors of
weight $\lambda$ in $V^{\otimes k}$.  Then any endomorphism of
$M^\lambda$ may be identified with an element of $U(\fq(V))$; the
restriction of this endomorphism on $R^\lambda$ may be identified with
an element of $U(\fh)$.

By Theorem 3.3.1 the endomorphism algebra of $R^\lambda$ is generated
by the $F_{ii}$ and to prove the corollary, it suffices to verify it
for these elements. We have
$$
F_{ii}(\kappa_{t}v_{t})= \kappa_{t}(F_{ii}v_{t})=
\kappa_{t}(p_{t}^iv_{t})=p_{t}^i\kappa_{t}(v_{t})
$$
\end{proof}

\ssbegin{3.3.3}{Corollary} Let $\varphi: M^\lambda\tto M^\lambda$ be
an $H_{k}$-module endomorphism given by the formula
$\varphi(v_{t})=\rho_{t}\kappa_{t}(v_{t})$.  Then $\varphi
|_{R^{\lambda}}=c\in \Cee$.
\end{Corollary}

\begin{proof} Let us show that $\varphi$ commutes with the
endomorphisms $F_{ii}$. Indeed,
$$
\begin{array}{l}
\varphi \cdot F_{ii}(v_{t})= \varphi(p_{t}^iv_{t})= p_{t}^i\varphi(v_{t})=
p_{t}^i\rho_{t}\kappa_{t}(v_{t});\\
F_{ii}\cdot \varphi (v_{t})= F_{ii}\rho_{t}\kappa_{t}(v_{t})=
\rho_{t}\kappa_{t}(F_{ii}v_{t})=\\
\rho_{t}\kappa_{t}(p_{t}^iv_{t})=\rho_{t}p_{t}^i\kappa_{t}(v_{t})=
p_{t}^i\rho_{t}\kappa_{t}(v_{t}).
\end{array}
$$
The latter identity holds thanks to the fact that $p_{t}^i$ commutes
with $\rho_{t}$.

Thus, $\varphi \cdot F_{ii}(v_{t})= F_{ii}\cdot \varphi (v_{t})$ and,
since the elements $v_{t}$ generate the $H_{k}$-module $M^\lambda$, we
have
$$
\varphi \cdot F_{ii}= F_{ii}\cdot \varphi .
$$
Since $\varphi(R^\lambda)\subset R^\lambda$ for any endomorphism
$\varphi $ of $M^\lambda$, it follows that $\varphi|_{R^\lambda}$ is
an element from the centralizer of $R^\lambda$.  But by Theorem 3.3.1
the centralizer is the Clifford algebra $C_{l(\lambda)}$.  But
$\varphi$ is an even central element of $C_{l(\lambda)}$, hence,
$\varphi $ is a constant.
\end{proof}

\ssbegin{3.3.4}{Corollary} Set $e_{t}=\kappa_{t}\rho_{t}$. Then
$$
e_{t}^2=c\cdot e_{t}\text{ for } c\in\Cee, c\neq 0
$$
and the algebra $e_{t}H_{k}e_{t}$ is isomorphic to $C_{l(\lambda)}$
and is generated by the $p_{t}^i$ for $1\leq i\leq l(\lambda)$.
\end{Corollary}

\begin{proof} Thanks to Corollary 3.3.3
$$
\varphi (\kappa_{t}(v_{t}))=\kappa_{t}\rho_{t}\kappa_{t}(v_{t})=\\
c\kappa_{t}(v_{t})
$$
or, equivalently,
$$
e_{t}^2=c\cdot e_{t}.
$$
Since the constant term of $\kappa_{t}\rho_{t}$, equal to the
coefficient of $v_{t}$ in $\kappa_{t}(v_{t})$, is nonzero, as follows
from the proof of Theorem 3.3.1, then the routine arguments with the
help of the bases (cf.  \cite{W}) shows that $c\neq 0$.  \end{proof}

Moreover, the algebra $e_{t}H_{k}e_{t}$ is {\it anti}-isomorphic to
the algebra of $H_{k}$-endomorphism of the submodule of $M^\lambda$
generated by $\kappa_{t}(v_{t})$.  Denote the latter module by
$R^\lambda_{1}$.  Thanks to Corollary 3.3.2 if $\varphi$ is an
endomorphism of $R^\lambda$, then $\varphi(R^\lambda_{1})\subset
R^\lambda_{1}$ and the restriction map $\varphi|_{\kappa_{t}(v_{t})}$
determines an antiisomorphism of $\End_{H_{k}}(R^\lambda)$ into the
algebra generated by the $p_{t}^i$ for $1\leq i\leq l(\lambda)$.  But
the latter algebra is isomorphic to $C_{(\lambda)}$; hence,
$\varphi|_{\kappa_{t}(v_{t})}$ is an anti-isomorphism and
$R^\lambda_{1}=R^\lambda$.

\section*{\S 4. The Specht modules over $\fA_{k}$}
In this section we construct analogs of the modules  $M^\lambda$ and
the Specht modules $R^\lambda$ for $\fA_{k}$.

First, we need the following statement.

\ssbegin{4.1}{Lemma} Let $\pi_{1}=\tau_{12}$, \dots ,
$\pi_{2}=\tau_{13}+\tau_{23}$, \dots , $\pi_{k}=\sum_{\alpha
<k}\tau_{\alpha k}$ be odd analogs of Jucys--Murphy's elements. Then
$$
e_{k}=\mathop{\prod}\limits_{i\geq 2}\frac{2}{i(i-1)}\pi_{i}^2
$$
is an idempotent and $e_{k}\fA_{k}e_{k}$ is isomorphic to the
Clifford algebra with $k-1$ generators.
\end{Lemma}

\begin{proof} It is easy to verify by induction that
$$
e_{k}=\frac{1}{k!}\sum_{\alpha\geq
0}(-1)^{\alpha}2^{k-\alpha-1}\Sigma_{2\alpha+1},
$$
where $\Sigma_{2\alpha+1}$ is the sum of all elements from $\fA_{k}$
of the form $\tau_{i_{1}i_{2}}\tau_{i_{2}i_{3}}\dots
\tau_{i_{2\alpha}i_{2\alpha+1}}$.  This implies that $e_{k}$ is a
central element that does not vary under the replacement of the
sequence $\pi_{i}=\sum_{\alpha <i}\tau_{\alpha i}$ with
$\pi_{\sigma(i)}=\sum_{\alpha <i}\tau_{\sigma(\alpha) \sigma(i)}$ for
any $i$ and $\sigma\in\fS_{k}*$.  It is not difficult to verify that
$$
(\tau_{12}+\tau_{23}+\tau_{31})\cdot(2-\tau_{12}\tau_{23}+
\tau_{13}\tau_{32})=0.
$$
Since
$$
2-\tau_{12}\tau_{23}+\tau_{13}\tau_{32}=\pi_{2}^2=
(\tau_{13}+\tau_{23})^{2}
$$
is a factor in the expression for $e_{k}$, it follows that
$$
(\tau_{12}+\tau_{23}+\tau_{31})e_{k}=0.
$$
From symmetry considerations
$$
(\tau_{ij}+\tau_{jl}+\tau_{li})e_{k}=0\text{ for any distinct } i, j,
l\in\{1, \dots k\}.
$$
Further on,
$$
\begin{array}{l}
\pi_{i}^2= (\sum\tau_{\alpha i})^{2}=i-1-\sum(\tau_{\alpha
\beta}\tau_{\beta i}+\tau_{\beta\alpha }\tau_{\alpha i})=\\
i-1-\left (\sum(1+\tau_{\alpha
\beta}\tau_{\beta i}+\tau_{\beta\alpha }\tau_{\alpha i})-1\right )=\\
i-1+\frac{1}{2}(i-1)(i-2)-
\sum\tau_{\alpha \beta}(\tau_{\alpha \beta}+\tau_{\beta i}+
\tau_{i\alpha }).\\
\end{array}
$$
Hence, $\pi_{i}^2e_{k}=\frac{1}{2}(i-1)(i-2)e_{k}$ and, therefore,
$$
e_{k}^2=\left(\mathop{\prod}\limits_{i\geq 2}\frac{2}{i(i-1)}
\pi_{i}^2\right )e_{k}= e_{k}.
$$

Furthermore, since $e_{k}$ is a central element, then
$e_{k}\fA_{k}e_{k}=\fA_{k}e_{k}$. Let $I$ be the ideal in $\fA_{k}$
generated by the elements $\tau_{ij}+\tau_{jl}+\tau_{li}$; then
$e_{n}I=0$. Let $\bar \fA= \fA_{k}/I$. In $\bar \fA$, then, the
following relations hold:
$$
\tau_{12}=\pi_{2}, \; \tau_{23}=\frac{1}{2}(\pi_{3}-\pi_{2}), \;
\tau_{34}=\frac{1}{3}(\pi_{4}-\pi_{3}), \; \dots ,\; \tau_{k-1,
k}=\frac{1}{k-1}(\pi_{k}-\pi_{k-1}).
$$
Hence,  $\bar \fA$ is generated by the $\pi_{i}$. Above we showed that
$\pi_{i}^2=\frac{1}{2}(i-1)(i-2)$ in $\bar \fA$. So $\bar \fA$ is the
Clifford algebra genereated by the images of the $\pi_{i}$ for $2\leq
i\leq k$.

Further on, $(1-e_{k})I=I$; so $\fA_{k}(1-e_{k})\supset I$, but since
$\bar \fA$ is a Clifford algebra (and, in particular, is simple), then
$\fA_{k}(1-e_{k})=I$ and, therefore, $\fA_{k}e_{k}\cong
\fA_{k}/(1-e_{k})\fA_{k}\cong \fA_{k}/I\cong \bar \fA$, the
Clifford algebra with $k-1$ generators.
\end{proof}

Let $\lambda$ be a strict partition and $t$ a $\lambda$-tableau. For
$1\leq i\leq l(\lambda)$ set
$$
p_{t}^{i}= \sum p_{\alpha}, \text{ where $\alpha$ runs over the
entries of the $i$-th row of $t$}
$$
and let $\sigma_{i}$ be the element from $\fA_{k}$ constructed as in
Lemma 4.1 from the sequence of numbers that stand in the $i$-th row.

Let $V$ be an $(n, n)$-dimensional superspace with $n\geq l(\lambda)$
and $W=V^{\otimes k}$, where $k=\sum \lambda_{i}$. Let $M^t$ be the
subspace of $W$ spanned by the vectors $w\in W$ such that
$$
E_{ii}w= \lambda_{i} w; \quad
F_{ii}w= p_{t}^{i} w\text{ for  $1\leq i\leq l(\lambda)$}.
$$
Set also $\sigma_{t}=\mathop{\prod}\limits_{1\leq i\leq
l(\lambda)}\sigma_{i}$; let
$R^t$ be the set of highest weight vectors that belong to $M^t$.

\ssbegin{4.2}{Theorem} {\em a)} As $\fA_{k}$-module, $M^t$ is
isomorphic to $\fA_{k}\sigma_{t}$.

{\em b)} The $\fA_{k}$-submodule in $M^t$ generated by the
$\kappa_{t'}v_{t}$, where $t'$ has the same rows as $t$, is an
isotypical one and its centralizer is isomorphic to the Clifford
algebra $C_{k-l(\lambda)}$.
\end{Theorem}

\begin{proof} Since $\sigma_{t}v_{t}=c\cdot v_{t}$, heading a) is
clear.

The facts $v_{t}\in M^t$ and $\fA_{k}v_{t}\subset R^t$ are also clear.

Let us show first of all that the centralizer of $R^t$ is isomorphic
to the Clifford algebra $C_{k-l(\lambda)}$.

Obviously, $\fA_{k}$ and $C_{k}\otimes U(\fq(V))$ form a Howe-dual
pair in $W$.  By Theorem 3.3.1 the centralizer of the $H_{k}$-module
$R^\lambda$ is the Clifford algebra generated by the $F_{ii}$ for
$1\leq i\leq l(\lambda)$.  Therefore, the centralizer of the
$\fA_{k}$-module $R^\lambda$ is the Clifford algebra
$C_{k+l(\lambda)}$ generated by the $F_{ii}$ for $1\leq i\leq
l(\lambda)$ and $p_{1}$, \dots , $p_{k}$.

The condition $F_{ii}w=p_{t}^{i}w$ is equivalent to another one:
$ew=w$ for
$$
e=\mathop{\prod}\limits_{1\leq i\leq l(\lambda)}\frac{1}{2}\left (1-
\frac{1}{\lambda_{i}}p_{t}^{i}F_{ii}\right )\text{ and } e^{2}=e.
$$
Therefore, the centralizer of $R^t$ is isomorphic to
$eC_{k+l(\lambda)}e$.

Let $C_{k-l(\lambda)}$ be the subalgebra of
$C_{k}$ generated by the $p_{i}-p_{j}$ for $i$, $j$ that belong to
the same row of $t$. Then it is not difficult to see that
$eC_{k+l(\lambda)}e\cong C_{k-l(\lambda)}$ and this proves that the
centralizer of $R^t$ is isomorphic to $C_{k-l(\lambda)}$.

Let us prove now that the submodule $R^t_{1}$ of $R^t$ generated by
the elements $\kappa_{t'}v_{t}$, where $t'$ has, up to
permutations, the same rows as $t$, is isomorphic to $R^t$. To this
end it suffices to prove that any endomorphism of $R^t$ sends $R^t_{1}$
into itself.

But, indeed, any endomorphism of $R^t$ is a multiplication by
$p_{i}-p_{j}$ for $i$, $j$ that belong to the same row of $t$.  Now,
it suffices to verify that $(p_{i}-p_{j})\kappa_{t'}v_{t}\in
R^t_{1}$. But
$$
(p_{i}-p_{j})\kappa_{t'}v_{t}=
\sqrt{2}\tau_{ij}s_{ij}\kappa_{t'}v_{t}=
-\sqrt{2}\tau_{ij}\kappa_{s_{ij}(t')}v_{t'}.
$$
But the rows of $s_{ij}t'$ consist of the same elements that
constitute the rows of $t$.
\end{proof}

\end{document}